\pdfoutput=1
\documentclass[11pt]{article}
\usepackage[left=1in,right=1in,top=1in,bottom=1in]{geometry}
\usepackage{times}
\usepackage{expl3}
\usepackage{cite}
\usepackage[table]{xcolor}
\usepackage{multirow}
\usepackage{stackengine} 
\usepackage{hhline}
\usepackage{lipsum}
\usepackage{titlesec}
\usepackage{wrapfig}
\usepackage{enumerate}
\usepackage{epsfig}
\usepackage{amsmath}
\usepackage{tabularx}
\usepackage{array}
\usepackage{booktabs}
\usepackage{enumitem}
\usepackage{bbm}
\usepackage{calc}
\usepackage{graphicx}
\usepackage{amsmath}
\usepackage[title]{appendix}
\usepackage{amssymb}
\usepackage{epstopdf}
\usepackage{boldline}
\usepackage{arydshln}
\usepackage{calligra}
\usepackage{bm}
\usepackage{url}
\usepackage{blindtext}
\usepackage{accents}

\newcommand{\define}{\stackrel{\mbox{\tiny def}}{=}}

\newtheorem{definition}{Definition}
\newtheorem{theorem}{Theorem}
\newtheorem{corollary}{Corollary}

\newtheorem{example}{Example}
\newtheorem{remark}{Remark}
\usepackage{mathtools}
\usepackage{epstopdf}
\usepackage{balance}
\usepackage{thmtools}
\usepackage{thm-restate}
\usepackage{hyperref}
\usepackage{cleveref}
\usepackage[mathscr]{euscript}

\usepackage[ruled,vlined]{algorithm2e}
\include{pythonlisting}
\newcommand{\ostar}{\mathbin{\mathpalette\make@circled\star}}

\makeatletter
\newcommand{\removelatexerror}{\let\@latex@error\@gobble}
\makeatother
\makeatletter
\newcommand*{\rom}[1]{\expandafter\@slowromancap\romannumeral #1@}
\makeatother

\ExplSyntaxOn
\newcommand\latinabbrev[1]{
  \peek_meaning:NTF . {% Same as \@ifnextchar
    #1\@}%
  { \peek_catcode:NTF a {% Check whether next char has same catcode as \'a, i.e., is a letter
      #1.\@ }%
    {#1.\@}}}
\ExplSyntaxOff

\titleclass{\subsubsubsection}{straight}[\subsubsection]

\begin{document}
\vspace{1cm}
\title{Chernoff Bounds for Tensor Expanders on Riemannian Manifolds Using Graph Laplacian Approximation}\vspace{1.8cm}
\author{Shih~Yu~Chang
% <-this % stops a space
\thanks{Shih Yu Chang is with the Department of Applied Data Science,
San Jose State University, San Jose, CA, U. S. A. (e-mail: {\tt
shihyu.chang@sjsu.edu}). 
           }}

\maketitle

\begin{abstract}
This paper addresses the advancement of probability tail bound analysis, a crucial statistical tool for assessing the probability of large deviations of random variables from their expected values. Traditional tail bounds, such as Markov's, Chebyshev's, and Chernoff bounds, have proven valuable across numerous scientific and engineering fields. However, as data complexity grows, there is a pressing need to extend tail bound estimation from scalar variables to high-dimensional random objects. Existing studies often rely on the assumption of independence among high-dimensional random objects, an assumption that may not always be valid. Building on the work of researchers like Garg et al. and Chang, who employed random walks to model high-dimensional ensembles, this study introduces a more generalized approach by exploring random walks over manifolds. To address the challenges of constructing an appropriate underlying graph for a manifold, we propose a novel method that enhances random walks on graphs approximating the manifold. This approach ensures spectral similarity between the original manifold and the approximated graph, including matching eigenvalues, eigenvectors, and eigenfunctions. Leveraging graph approximation technique proposed by Burago et al. for manifolds, we derive the tensor Chernoff bound and establish its range for random walks on a Riemannian manifold according to the underlying manifold's spectral characteristics.
\end{abstract}

\begin{keywords}
Tensors, random walk, tail bound, Chernoff bound, manifold, graph.
\end{keywords}

\section{Introduction}\label{sec: Introduction}

Probability tail bound analysis is a powerful statistical tool used to quantify the probability that a random variable deviates significantly from its expected value. Tail bounds provide inequalities, such as Markov's, Chebyshev's, and Chernoff bounds, which help estimate the likelihood of extreme events in the distribution's tails~\cite{ledoux2001concentration}. These bounds are crucial in various fields of science and engineering for ensuring reliability and safety~\cite{lopes2022sharp}. In science applications, tail bounds are used in experimental data analysis to identify outliers and assess risks. For example, in physics, they help evaluate the likelihood of rare events in quantum systems. In engineering aspect, tail bounds are applied in fields such as telecommunications to ensure that signal loss or delays remain within acceptable limits~\cite{li2021ultra}. Tail bounds also play a significant role in machine learning and computer science for algorithm performance analysis, ensuring that computational processes remain efficient and reliable under unpredictable conditions~\cite{tropp2015introduction}.

Enhancing tail bound estimation theory from scalar to high-dimensional random objects is vital due to the increasing complexity of data and problems in modern applications. Fields such as machine learning, finance, genomics, and physics involve numerous interdependent variables, which cannot be effectively captured by scalar random variables. High-dimensional tail bounds analysis equips us with the necessary tools to model, analyze, and comprehend the interactions among multiple variables. Moreover, real-world data often exist in multi-dimensional forms, including images, videos, and multi-sensor measurements~\cite{chang2022TWF,chang2023TLS,chang2023TEKF}. Extending tail bound estimation theory to these dimensions allows us to manage the complexities of such data, enabling more accurate predictions, inferences, and decisions. Applications of high-dimensional tail bound analysis in vector, matrix, and tensor formats are evident in control system design, optimization theory, randomized algorithm design, numerical analysis, and quantum mechanics~\cite{louart2018random,dulian2023matrix}.

Many existing analyses of tail bounds for high-dimensional objects assume that the ensemble of random objects is independent~\cite{chang2022convenient,chang2021TensorExp}. However, this independence assumption may not be valid in all cases. An alternative approach involves developing tail bounds under non-independence assumptions for random variables. Gillman~\cite{gillman1998chernoff}, along with subsequent refinements~\cite{chung2012chernoff, rao2017sharp}, replaced the independence assumption with Markov dependence. Their work can be summarized as follows: Let $\mathscr{G}$ be a regular $\lambda$-expander graph with a vertex set $\mathscr{V}$, and let $g: \mathscr{V} \rightarrow \mathbb{C}$ be a bounded function. Consider a stationary random walk $v_1, v_2, \ldots, v_{K}$ of length $K$ on $\mathscr{G}$. They demonstrated that:
\begin{eqnarray}\label{eq:Chernoff bound Markov rvs}
\mathrm{Pr} \left( \left\vert \frac{1}{K} \sum\limits_{j=1}^{K} g(v_i) - \mathbb{E}[g] \right\vert \geq \vartheta \right) \leq 2 \exp( - \Omega (1 - \lambda) K \vartheta^2),
\end{eqnarray}
where the parameter $\lambda$ corresponds to the second-largest eigenvalue of the transition matrix of the underlying graph $\mathfrak{G}$. The bound presented in Eq.~\eqref{eq:Chernoff bound Markov rvs} is known as the "Expander Chernoff Bound." This naturally leads to the extension of Eq.~\eqref{eq:Chernoff bound Markov rvs} into the "Matrix Expander Chernoff Bound." Wigderson and Xiao initially attempted to establish partial results for the "Matrix Expander Chernoff Bound" in~\cite{wigderson2008derandomizing}, with Garg et al.~\cite{garg2018matrix} later providing a comprehensive solution. In~\cite{chang2021TensorExp}, the author extends the matrix expander Chernoff bound by Garg et al.~\cite{garg2018matrix} to tensor expander Chernoff bounds, utilizing new tensor norm inequalities based on log-majorization techniques. Theose bounds derived in~\cite{chang2021TensorExp} enable (1) the expansion from matrices to tensors, (2) the application to any polynomial function of the summed random objects, (3) the use of the Ky Fan norm beyond solely eigenvalues, and (4) the elimination of the zero-sum restriction in mapped random walks to tensors.

Considering random walks over manifolds is crucial for understanding complex systems where data naturally reside in curved spaces rather than flat, Euclidean spaces. Manifolds provide a mathematically rigorous framework to study phenomena such as diffusion, heat distribution, and biological processes, offering a richer context than traditional Euclidean settings~\cite{lee2018introduction}. This approach is particularly important in fields like machine learning, where data often exists on nonlinear structures, and in physics, where the geometry of space itself is curved. Analyzing random walks on manifolds helps develop algorithms and models that are sensitive to the intrinsic geometry of data, leading to more accurate and effective solutions in manifold learning~\cite{izenman2012introduction}.

The exploration of random walk behavior on manifolds dates back to the 1940s. We will discuss several key contributions here. In~\cite{roberts1960random}, the authors first examined random walks on spheres by conceptualizing them as sequences of random steps with directions uniformly distributed from the origin. Each step's length was either fixed or governed by a specified probability distribution, with the lengths allowed to vary during the walk. This research provided a straightforward solution for determining the endpoint distribution of any random walk, leveraging the commutative nature and shared eigenfunctions of the individual steps. The study was extended to random walks on general Riemannian manifolds, noting that while commutativity is not universally applicable, it is present in completely harmonic spaces and certain others, aligning the methodology with that of the sphere in these instances. 

In~\cite{jorgensen1975central}, the author investigated random walks on general Riemannian manifolds and explored the limiting behavior of sequences of such walks. The study showed that under reasonable conditions, these sequences converge to a diffusion process on the manifold, with Brownian motion processes emerging as the limits of random walks with identically distributed steps through the application of semigroup theory. The approach employed semigroup methods akin to those in earlier works~\cite{hunt1956semi}. The semigroup method involved defining sub-probability measures on the tangent spaces of the manifold, constructing a random walk, and subsequently an associated Markov process. By imposing certain conditions on these sequences, the author in~\cite{jorgensen1975central} demonstrated that the corresponding Markov processes converge weakly to a diffusion process.

In~\cite{varopoulos1984brownian}, the author explored the relationship between canonical Brownian motion on a manifold and the admissible random walks on its grid approximation. The paper established that the recurrence of Brownian motion on the manifold is equivalent to the recurrence of the defined random walks on its grid approximation, provided the approximation is sufficiently close for any point of the manifold to its nearest grid point. More recent studies on random walk behaviors over Riemannian manifolds can also be found in~\cite{mangoubi2018rapid,herrmann2022efficient,ma2021geodesic}. However, these works do not offer spectral information about the transition matrix that characterizes random walk behavior over the underlying manifold. From Eq.~\eqref{eq:Chernoff bound Markov rvs} (scalar case) and the derivation of high-dimensional Chernoff bounds (matrix~\cite{garg2018matrix} and tensor~\cite{chang2021TensorExp}), it is evident that spectral information, particularly the second largest eigenvalue of the transition matrix, is crucial for determining the Chernoff bound for various types of random objects.

Constructing an appropriate underlying graph for random walks over a manifold with desired graph spectrum properties is challenging. To address this, we adopt a reverse approach, enhancing random walks on a graph to closely approximate the manifold by ensuring spectral closeness, encompassing eigenvalues, eigenvectors, and eigenfunctions, between the original manifold and the approximated random walk graph. Utilizing the methodology described in~\cite{burago2015graph}, we construct a graph $\mathscr{G}_{\mathscr{M}}$ that approximates the manifold $\mathscr{M}$ such that the spectral difference between the Laplacian matrix of $\mathscr{G}_{\mathscr{M}}$ and the Laplace-Beltrami operator of $\mathscr{M}$ remains within a controllable bound dictated by the manifold's properties. This allows us to derive the tensor Chernoff bound and its range for random walks over a Riemannian manifold based on the spectral properties of the manifold $\mathscr{M}$.

The remainder of this paper is structured as follows: In Section~\ref{sec: Graph Approximation of a Riemannian Manifold}, we apply weighted graph approximation theory to Riemannian manifolds to estimate the spectral bounds of transition matrices for random walks on Riemannian manifolds. In Section~\ref{sec: Tensor Expander Chernoff Bounds over a Riemannian Manifold}, we establish the tail bounds for tensor expander Chernoff bounds for random walks over a Riemannian manifold.

\noindent \textbf{Nomenclature:} To simplify notation, let $\mathbb{I}_1^M$ be defined as $\prod_{i=1}^M I_i$, where $I_i$ indicates the size of the $i$-th dimension of a tensor. The term \emph{Hermitian tensor} is specified by Definition~3 in~\cite{chang2022convenient}. The symbols $\lambda_{\max}(\mathcal{H})$ denote the largest eigenvalue of a Hermitian tensor $\mathcal{H}$, as explained in~\cite{chang2022convenient}. The Ky Fan-like $k$-norm of the Hermitian tensor $\mathcal{H}$ is represented as $\left\Vert \mathcal{H}\right\Vert_{(k)}$~\cite{chang2021TensorExp}.

\section{Graph Approximation of a Riemannian Manifold}\label{sec: Graph Approximation of a Riemannian Manifold}

In this section, we will utilize the weighted graph approximation theory for Riemannian manifolds, as proposed by~\cite{burago2015graph}, to estimate the lower and upper bounds of the spectrum of transition matrices for random walks on manifolds.

We will begin with Definition~\ref{def: wighted graph to app manifold M} of a weighted graph used to approximate a Riemannian manifold.  
\begin{definition}\label{def: wighted graph to app manifold M}
Gievn a Riemannian manifold $(\mathscr{M}, g)$ with dimension $n$ consists of a smooth manifold $\mathscr{M}$ equipped with a Riemannian metric $g$. The metric $g$ is a smoothly varying positive-definite symmetric bilinear form on the tangent space of $\mathscr{M}$ at each point. We will establish a weighted graph, denoted by $\mathscr{G}_{\mathscr{M}}(\epsilon,\mu,\kappa)=(\mathscr{V}, \mathscr{E}, \mathscr{W})$, where $\mathscr{V}={v_i}$ is the sef of vertices sampled from the manifold $\mathscr{M}$ for $i=1,2,\ldots,N$ and $\mathscr{E}=\{e_{i,j}=(v_i, v_j)\}$ is the set of edges, to approximate the manifold $\mathscr{M}$ such that 
\begin{itemize}
\item All balls with centers $v_i$ with radius $\epsilon$, represented by $B_{\epsilon}(v_i)$, can cover $\mathscr{M}$, i.e., $\mathscr{M} \subset \bigcup\limits_{i=1}^N B_{\epsilon}(v_i)$;   
\item The measure $\mu$ on the set $\{v_i\}$ will assign the measure $\mu_i$ to the volume of the space $V_i$, where $V_i$ satisfies $\mathscr{M}=\bigcup\limits_{i=1}^N V_i$ and $V_i \subset B_\epsilon(v_i)$; 
\item The edge $e_{i,j}=(v_i, v_j)$ is formed if $d_g(v_i, v_j) < \kappa$, and the weight $w_{i,j} \subset \mathscr{W}$ for the edge $e_{i,j}\in \mathscr{E}$ is determined by 
\begin{eqnarray}\label{eq: edge weight def}
w_{i,j}\define\frac{2(n+2)\Gamma(1+n/2)}{\pi^{n/2}\kappa^{n+2}}\mu_i \mu_j,
\end{eqnarray}
where $\Gamma$ is the Gamma function. 
\end{itemize}
\end{definition}

The weighted graph $\mathscr{G}_{\mathscr{M}}(\epsilon,\mu,\kappa)=(\mathscr{V}, \mathscr{E}, \mathscr{W})$ always exists because the vertex set $\mathscr{V}$ can be derived from the centers $v_i$ of the Voronoi decomposition cells $V_i$ of a manifold $\mathscr{M}$, where each Voronoi cell $V_i \subseteq B_{\epsilon}(v_i)$. 

From the weighted graph $\mathscr{G}_{\mathscr{M}}(\epsilon,\mu,\kappa)=(\mathscr{V}, \mathscr{E}, \mathscr{W})$ provided by Definition~\ref{def: wighted graph to app manifold M}, we can form a diagonal matrix based on $\mathscr{G}_{\mathscr{M}}(\epsilon,\mu,\kappa)$ to represent the sum of the weights of the edges connected to vertex $v_i$. This diagonal matrix, a.k.a, degree matrix, is represented by $\bm{D}_{\mathscr{G}_{\mathscr{M}}}$, which can be expressed by
\begin{eqnarray}\label{eq: degree matrix def}
\bm{D}_{\mathscr{G}_{\mathscr{M}}}&=&[d_{i,i}]=\left[\sum\limits_{j=1}^N w_{i,j}\right], 
\end{eqnarray}
where $i=1,2,\ldots,N$. We also can form an adjancy matrix based on $\mathscr{G}_{\mathscr{M}}(\epsilon,\mu,\kappa)$ to represent the edge bewteen vertices $v_i$ and $v_j$. The adjancy matrix is denoted by $\bm{A}_{\mathscr{G}_{\mathscr{M}}}$, which can be expressed by
\begin{eqnarray}\label{eq: adjancy matrix def}
\bm{A}_{\mathscr{G}_{\mathscr{M}}}&=&[a_{i,j}]=w_{i,j}, 
\end{eqnarray}
where $i, j =1,2,\ldots,N$ but $i \neq j$. From Eq.~\eqref{eq: degree matrix def} and Eq.~\eqref{eq: adjancy matrix def}, we can form the Laplacian matrix for the weighted graph $\mathscr{G}_{\mathscr{M}}(\epsilon,\mu,\kappa)$ by
\begin{eqnarray}\label{eq: Laplacian for app weighted graph}
\bm{L}_{\mathscr{G}_{\mathscr{M}}}&=&\bm{D}_{\mathscr{G}_{\mathscr{M}}} - \bm{A}_{\mathscr{G}_{\mathscr{M}}}.
\end{eqnarray}

By multiplying $\bm{D}^{-1}_{\mathscr{G}_{\mathscr{M}}}$ at both sides of Eq.~\eqref{eq: Laplacian for app weighted graph}, we will get 
\begin{eqnarray}\label{eq: Laplacian for app norm weighted graph}
\tilde{\bm{L}}_{\mathscr{G}_{\mathscr{M}}}&\define&\bm{D}^{-1}_{\mathscr{G}_{\mathscr{M}}}\bm{L}_{\mathscr{G}_{\mathscr{M}}}\nonumber \\
&=&\bm{I} - \bm{D}^{-1}_{\mathscr{G}_{\mathscr{M}}}\bm{A}_{\mathscr{G}_{\mathscr{M}}}\nonumber \\
&=&\bm{I} - \bm{P}_{\mathscr{G}_{\mathscr{M}}}, 
\end{eqnarray}
where we set $\bm{P}_{\mathscr{G}_{\mathscr{M}}}\define\bm{D}^{-1}_{\mathscr{G}_{\mathscr{M}}}\bm{A}_{\mathscr{G}_{\mathscr{M}}}$ at the last equality. The matrix $\bm{P}_{\mathscr{G}_{\mathscr{M}}}$ will be used as the transion matrix for random walks over vertices $\mathscr{V}$ of the weighted graph $\mathscr{G}_{\mathscr{M}}(\epsilon,\mu,\kappa)$, which is a discrete approximation of the underlying manifold $\mathscr{M}$. 

Let $\lambda_{\bm{L}_{\mathscr{G}_{\mathscr{M}}},i}$ be eigenvalues of the Laplacian matrix $\bm{L}_{\mathscr{G}_{\mathscr{M}}}$ and $\lambda_{\bm{P}_{\mathscr{G}_{\mathscr{M}}},i}$ be eigenvalues of the transition matrix $\bm{P}_{\mathscr{G}_{\mathscr{M}}}$, then, from Eq.~\eqref{eq: Laplacian for app norm weighted graph}, we have
\begin{eqnarray}\label{eq: lambda P g M and lambda L g m relation}
\lambda_{\bm{P}_{\mathscr{G}_{\mathscr{M}}},i}&=&1-\frac{\lambda_{\bm{L}_{\mathscr{G}_{\mathscr{M}}},i}}{\sum\limits_{j=1}^N w_{i,j}},
\end{eqnarray}
where $w_{i,j}$ is definied by Eq.~\eqref{eq: edge weight def}.

From Theorem 1 and its conditions in~\cite{burago2015graph}, we have 
\begin{eqnarray}\label{eq: spectrum bounds by graph discrete paper}
\left\vert \lambda_{\bm{L}_{\mathscr{G}_{\mathscr{M}}},i}-\lambda_{\mathscr{M},i}\right\vert
&\leq& C_{n,D_{\mathscr{M}},r_{\mathscr{M}}}\left[(\epsilon/\kappa + K_{\mathscr{M}}\kappa^2)\lambda_{\mathscr{M},i}+\kappa \lambda_{\mathscr{M},i}^{3/2}\right],
\end{eqnarray}
where $\lambda_{\mathscr{M},i}$ are eigenvalues of the Laplace-Beltrami operator on the manifold $\mathscr{M}$ and $C_{n,D_{\mathscr{M}},r_{\mathscr{M}}}$ is the constant for the $i$-th eigenvalue associated to the underlying manifold properties of $\mathscr{M}$, which are the diameter $D_{\mathscr{M}}$ and the injectivity radius $r_{\mathscr{M}}$. 

\section{Tensor Expander Chernoff Bounds over a Riemannian Manifold}\label{sec: Tensor Expander Chernoff Bounds over a Riemannian Manifold}

The main purpose of this section is to prove Theorem~\ref{thm:tensor expander}, which provides the tail bound for tensor expander Chernoff bounds over a Riemannian manifold. This Theorem~\ref{thm:tensor expander} is the extension of Theorem 4.4 in~\cite{chang2021TensorExp} by associating the random walks with the underlying manifold properties.

\begin{theorem}\label{thm:tensor expander}
Let $\mathscr{G}_{\mathscr{M}}(\epsilon,\mu,\kappa)=(\mathscr{V}, \mathscr{E}, \mathscr{W})$ be the approximation graph for $\mathscr{M}$ whose transition matrix has second largest eigenvalue $\lambda_{\bm{P}_{\mathscr{G}_{\mathscr{M}}},\tilde{i}}$ for some $\tilde{i} \in \{1,2,\ldots,N\}$, and let $g: \mathscr{V} \rightarrow \in \mathbb{C}^{I_1 \times \cdots \times I_M \times I_1 \times \cdots \times I_M}$ be a function. We assume following: 
\begin{enumerate}
\item For each $v_i \in \mathscr{V}$, $g(v_i)$ is a Hermitian tensor;
\item $\left\Vert g(v_i) \right\Vert \leq r$ for all $v_i \in \mathscr{V}$;
\item A nonnegative coefficients polynomial raised by the power $s \geq 1$ as $f: x \rightarrow (a_0 + a_1x +a_2 x^2 + \cdots +a_n x^n)^s$ satisfying $f \left(\exp \left( t   \sum\limits_{i=1}^{K} g(v_i) \right) \right) \succeq \exp \left( t f \left(  \sum\limits_{i=1}^{K} g(v_i) \right) \right) $ almost surely;
\item  For $\tau \in [\infty, \infty]$, we have constants $C$ and $\sigma$ such that $ \beta_0(\tau)  \leq \frac{C \exp( \frac{-\tau^2}{2 \sigma^2} ) }{\sigma \sqrt{2 \pi}}$, where $\beta_0(\tau)$ is the interpolation function used to establish tensor norm inequalities~\cite{chang2021TensorExp}.
\end{enumerate}
Then, we have 
\begin{eqnarray}\label{eq0:thm:tensor expander}
\mathrm{Pr}\left( \left\Vert f \left(  \sum\limits_{i=1}^{K} g(v_i) \right)  \right\Vert_{(k)} \geq\vartheta \right) \leq   \min\limits_{t > 0 } \left[ (n+1)^{(s-1)} e^{-\vartheta t} \left(a_0 K  +C \left( K + \sqrt{\frac{\mathbb{I}_1^M - K }{K}}\right)\cdot \right. \right. ~~~~~  \nonumber \\
\left. \left. \sum\limits_{l=1}^n a_l \exp( 8 K \overline{\lambda_{\bm{P}_{\mathscr{G}_{\mathscr{M}}},\tilde{i}}} + 2  (K +8 \overline{\lambda_{\bm{P}_{\mathscr{G}_{\mathscr{M}}},\tilde{i}}}) lsr t + 2(\sigma (K +8 \overline{\lambda_{\bm{P}_{\mathscr{G}_{\mathscr{M}}},\tilde{i}}}) lsr )^2 t^2  )  \right)\right]\nonumber \\
=   \min\limits_{t > 0 } \left[ (n+1)^{(s-1)} e^{-\vartheta t} \left(a_0 K  +C \left( K + \sqrt{\frac{\mathbb{I}_1^M - K }{K}}\right)\cdot \right. \right. ~~~~~~~~~~~~~~~~~~~~~~~~~~~~~~  \nonumber \\
\left. \left. \sum\limits_{l=1}^n a_l \exp\left( 8 K \frac{\lambda_{\bm{L}_{\mathscr{G}_{\mathscr{M}}},\tilde{i}}}{\sum\limits_{j=1}^N w_{\tilde{i},j}} + 2\left(K +8\frac{\lambda_{\bm{L}_{\mathscr{G}_{\mathscr{M}}},\tilde{i}}}{\sum\limits_{j=1}^N w_{\tilde{i},j}}\right) lsrt + 2\left(\sigma\left(K +8 \frac{\lambda_{\bm{L}_{\mathscr{G}_{\mathscr{M}}},\tilde{i}}}{\sum\limits_{j=1}^N w_{\tilde{i},j}}\right)lsr\right)^2 t^2  \right)\right)\right],
\end{eqnarray}
where $\overline{\lambda_{\bm{P}_{\mathscr{G}_{\mathscr{M}}},\tilde{i}}} = 1 - \lambda_{\bm{P}_{\mathscr{G}_{\mathscr{M}}},\tilde{i}}$. 
\end{theorem}
\textbf{Proof:}
By applying Eq.~\eqref{eq: lambda P g M and lambda L g m relation} to the eigenvalues of the transition matrix of the random walk over the manifold $\mathscr{M}$, as stated in Theorem 4.4 of~\cite{chang2021TensorExp}, the theorem is thereby proven.
$\hfill \Box$

From Eq.~\eqref{eq: spectrum bounds by graph discrete paper}, we have
\begin{eqnarray}\label{eq: spectrum U L bounds by graph discrete paper}
 \underbrace{\lambda_{\mathscr{M},i} - C_{n,D_{\mathscr{M}},r_{\mathscr{M}}}\left[(\epsilon/\kappa + K_{\mathscr{M}}\kappa^2)\lambda_{\mathscr{M},i}+\kappa \lambda_{\mathscr{M},i}^{3/2}\right]}_{\define \Phi_{\mathscr{L}}(\lambda_{\mathscr{M},i})}&\leq&
\lambda_{\bm{L}_{\mathscr{G}_{\mathscr{M}}},i} \nonumber \\
&\leq&  \underbrace{\lambda_{\mathscr{M},i} + C_{n,D_{\mathscr{M}},r_{\mathscr{M}}}\left[(\epsilon/\kappa + K_{\mathscr{M}}\kappa^2)\lambda_{\mathscr{M},i}+\kappa \lambda_{\mathscr{M},i}^{3/2}\right]}_{\define \Phi_{\mathscr{U}}(\lambda_{\mathscr{M},i})}.\nonumber \\
\end{eqnarray}
Thus, by applying the Extreme Value Theorem, we can establish both upper and lower bounds for the tail probability in Eq.~\eqref{eq0:thm:tensor expander}, as the eigenvalues \(\lambda_{\bm{L}_{\mathscr{G}_{\mathscr{M}}},i}\) can be bounded using Eq.~\eqref{eq: spectrum U L bounds by graph discrete paper}. Should additional constraints be imposed under Theorem~\ref{thm:tensor expander}, these bounds for the tail probability in Eq.~\eqref{eq0:thm:tensor expander} can be determined explicitly. Corollary~\ref{cor: tensor expander a i greater zero} will elaborate on these results.

\begin{corollary}\label{cor: tensor expander a i greater zero}
If all coefficients $a_i$ for the following mapping are nonnegative:
\begin{eqnarray}\label{eq1: cor: tensor expander a i greater zero}
f: x \rightarrow (a_0 + a_1x +a_2 x^2 + \cdots +a_n x^n)^s, 
\end{eqnarray}
and the following expression is monotone increasing with respect to $\lambda_{\bm{L}_{\mathscr{G}_{\mathscr{M}}},\tilde{i}}$ for $\lambda_{\bm{L}_{\mathscr{G}_{\mathscr{M}}},\tilde{i}}$ in the range provided by Eq.~\eqref{eq: spectrum U L bounds by graph discrete paper}: 
\begin{eqnarray}
\exp\left( 8 K \frac{\lambda_{\bm{L}_{\mathscr{G}_{\mathscr{M}}},\tilde{i}}}{\sum\limits_{j=1}^N w_{\tilde{i},j}} + 2\left(K +8\frac{\lambda_{\bm{L}_{\mathscr{G}_{\mathscr{M}}},\tilde{i}}}{\sum\limits_{j=1}^N w_{\tilde{i},j}}\right) lsrt + 2\left(\sigma\left(K +8 \frac{\lambda_{\bm{L}_{\mathscr{G}_{\mathscr{M}}},\tilde{i}}}{\sum\limits_{j=1}^N w_{\tilde{i},j}}\right)lsr\right)^2 t^2  \right). 
\end{eqnarray}

Then, we have the upper bound for the tail probability in Eq.~\eqref{eq0:thm:tensor expander} as:
\begin{eqnarray}\label{eq2: cor: tensor expander a i greater zero}
\mathrm{Pr}\left( \left\Vert f \left(  \sum\limits_{i=1}^{K} g(v_i) \right)  \right\Vert_{(k)} \geq\vartheta \right) \leq    \min\limits_{t > 0 } \left[ (n+1)^{(s-1)} e^{-\vartheta t} \left(a_0 K  +C \left( K + \sqrt{\frac{\mathbb{I}_1^M - K }{K}}\right)\cdot \right. \right. ~~~~~  \nonumber \\
\left. \left. \sum\limits_{l=1}^n a_l \exp\left( 8 K \frac{\Phi_{\mathscr{U}}(\lambda_{\mathscr{M},i})}{\sum\limits_{j=1}^N w_{i,j}} + 2  \left(K +8\frac{\Phi_{\mathscr{U}}(\lambda_{\mathscr{M},i})}{\sum\limits_{j=1}^N w_{i,j}}\right) lsr t + 2\left(\sigma\left(K +8 \frac{\Phi_{\mathscr{U}}(\lambda_{\mathscr{M},i})}{\sum\limits_{j=1}^N w_{i,j}}\right) lsr\right)^2 t^2\right)  \right)\right].
\end{eqnarray}
Additionaly, we also have the lower bound for the tail probability in Eq.~\eqref{eq0:thm:tensor expander} as:
\begin{eqnarray}\label{eq3: cor: tensor expander a i greater zero}
\mathrm{Pr}\left( \left\Vert f \left(  \sum\limits_{i=1}^{K} g(v_i) \right)  \right\Vert_{(k)} \geq\vartheta \right) \leq    \min\limits_{t > 0 } \left[ (n+1)^{(s-1)} e^{-\vartheta t} \left(a_0 K  +C \left( K + \sqrt{\frac{\mathbb{I}_1^M - K }{K}}\right)\cdot \right. \right. ~~~~~  \nonumber \\
\left. \left. \sum\limits_{l=1}^n a_l \exp\left( 8 K \frac{\Phi_{\mathscr{L}}(\lambda_{\mathscr{M},i})}{\sum\limits_{j=1}^N w_{i,j}} + 2  \left(K +8\frac{\Phi_{\mathscr{L}}(\lambda_{\mathscr{M},i})}{\sum\limits_{j=1}^N w_{i,j}}\right) lsr t + 2\left(\sigma\left(K +8 \frac{\Phi_{\mathscr{L}}(\lambda_{\mathscr{M},i})}{\sum\limits_{j=1}^N w_{i,j}}\right) lsr\right)^2 t^2\right)  \right)\right].
\end{eqnarray}
\end{corollary}
\textbf{Proof:}
Because all coeffecient $a_i$ in Eq.~\eqref{eq1: cor: tensor expander a i greater zero} are nonnegative, the summation of monotone functions $\exp\left( 8 K \frac{\lambda_{\bm{L}_{\mathscr{G}_{\mathscr{M}}},\tilde{i}}}{\sum\limits_{j=1}^N w_{\tilde{i},j}} + 2\left(K +8\frac{\lambda_{\bm{L}_{\mathscr{G}_{\mathscr{M}}},\tilde{i}}}{\sum\limits_{j=1}^N w_{\tilde{i},j}}\right) lsrt + 2\left(\sigma\left(K +8 \frac{\lambda_{\bm{L}_{\mathscr{G}_{\mathscr{M}}},\tilde{i}}}{\sum\limits_{j=1}^N w_{\tilde{i},j}}\right)lsr\right)^2 t^2  \right)$ with respect to $\lambda_{\bm{L}_{\mathscr{G}_{\mathscr{M}}},\tilde{i}}$ is monotone.  

From the first inequality of Eq.~\eqref{eq: spectrum U L bounds by graph discrete paper} and Theorem~\ref{thm:tensor expander}, we have Eq.~\eqref{eq2: cor: tensor expander a i greater zero}. Similarly, from the second inequality Eq.~\eqref{eq: spectrum U L bounds by graph discrete paper} and Theorem~\ref{thm:tensor expander}, we have Eq.~\eqref{eq3: cor: tensor expander a i greater zero}
$\hfill\Box$

\begin{example}\label{exp: ten exp over sphere n}
We consider an $n$-dimensional sphere with radius $R$, represented by $S^n(R)$, as our manifold $\mathscr{M}$, and the metric tensor $g_{i,j}$ induced hyperspherical coordinates can be expressed by
\begin{eqnarray}
g_{i,j}&=&\begin{bmatrix}
1 & 0 & 0 & \ldots & 0\\
0 & \sin^2\theta_1 & 0 & \ldots & 0\\
0 & 0 & \sin^2\theta_1 \sin^2\theta_2 & \ldots & 0\\
\vdots & \vdots & \vdots & \ddots & \vdots\\
0 & 0 & 0 & \ldots & \prod\limits_{k=1}^{n-1}\sin^2\theta_k \\
\end{bmatrix}.
\end{eqnarray} 
Thus, the eigenvalues of the Laplace-Beltrami operator on the $S^n(R)$ are $\lambda_{S^n(R),i} = (i-1)(i + n - 2)$ for $i=1,2,\ldots,n$. The sectional curvature of $S^n(R)$ is $\frac{1}{R^2}$.  

From Eq.~\eqref{eq: spectrum U L bounds by graph discrete paper}, we have
\begin{eqnarray}\label{eq1: exp: ten exp over sphere n}
\Phi_{\mathscr{U}}(\lambda_{S^n(R),i})&=&(i-1)(i + n - 2) + C_{n,D_{S^n(R)},r_{S^n(R)}} \nonumber \\
&&\cdot \left[\left(\epsilon/\kappa + \frac{1}{R^2}\kappa^2\right)(i-1)(i + n - 2)+\kappa ((i-1)(i + n - 2))^{3/2}\right],
\end{eqnarray}
and
\begin{eqnarray}\label{eq2: exp: ten exp over sphere n}
\Phi_{\mathscr{L}}(\lambda_{S^n(R),i})&=&(i-1)(i + n - 2)  - C_{n,D_{S^n(R)},r_{S^n(R)}} \nonumber \\
&&\cdot \left[\left(\epsilon/\kappa + \frac{1}{R^2}\kappa^2\right)(i-1)(i + n - 2)+\kappa ((i-1)(i + n - 2))^{3/2}\right]. 
\end{eqnarray}
If those conditions required by Corollar~\ref{cor: tensor expander a i greater zero} are satsified, from Corollar~\ref{cor: tensor expander a i greater zero}, we have the upper bound for the tail probability in Eq.~\eqref{eq0:thm:tensor expander} for $S^n(R)$ as:
\begin{eqnarray}\label{exp1: cor: tensor expander a i greater zero}
\mathrm{Pr}\left( \left\Vert f \left(  \sum\limits_{i=1}^{K} g(v_i) \right)  \right\Vert_{(k)} \geq\vartheta \right) \leq    \min\limits_{t > 0 } \Bigg[ (n+1)^{(s-1)} e^{-\vartheta t} \Bigg(a_0 K  +C \left( K + \sqrt{\frac{\mathbb{I}_1^M - K }{K}}\right)~~~~~~~~~~    \nonumber \\
\cdot \sum\limits_{l=1}^n a_l \exp\Bigg( 8 K \frac{\Phi_{\mathscr{U}}(\lambda_{S^n(R),\tilde{i}})}{\sum\limits_{j=1}^N w_{\tilde{i},j}} + 2  \left(K +8\frac{\Phi_{\mathscr{U}}(\lambda_{S^n(R),\tilde{i}})}{\sum\limits_{j=1}^N w_{\tilde{i},j}}\right) lsrt \nonumber \\
+ 2\left(\sigma\left(K +8 \frac{\Phi_{\mathscr{U}}(\lambda_{S^n(R),\tilde{i}})}{\sum\limits_{j=1}^N w_{\tilde{i},j}}\right) lsr\right)^2 t^2\Bigg)  \Bigg)\Bigg].~~~~~~~~~~~~~~~~~~~~~~~~~~
\end{eqnarray}
Additionaly, we also have the lower bound for the tail probability in Eq.~\eqref{eq0:thm:tensor expander} for $S^n(R)$ as:
\begin{eqnarray}\label{exp2: cor: tensor expander a i greater zero}
\mathrm{Pr}\left( \left\Vert f \left(  \sum\limits_{i=1}^{K} g(v_i) \right)  \right\Vert_{(k)} \geq\vartheta \right) \leq    \min\limits_{t > 0 } \Bigg[ (n+1)^{(s-1)} e^{-\vartheta t} \Bigg(a_0 K  +C \left( K + \sqrt{\frac{\mathbb{I}_1^M - K }{K}}\right)~~~~~~~~~~    \nonumber \\
\cdot \sum\limits_{l=1}^n a_l \exp\Bigg( 8 K \frac{\Phi_{\mathscr{L}}(\lambda_{S^n(R),\tilde{i}})}{\sum\limits_{j=1}^N w_{\tilde{i},j}} + 2  \left(K +8\frac{\Phi_{\mathscr{L}}(\lambda_{S^n(R),\tilde{i}})}{\sum\limits_{j=1}^N w_{\tilde{i},j}}\right) lsrt \nonumber \\
+ 2\left(\sigma\left(K +8 \frac{\Phi_{\mathscr{L}}(\lambda_{S^n(R),\tilde{i}})}{\sum\limits_{j=1}^N w_{\tilde{i},j}}\right) lsr\right)^2 t^2\Bigg)  \Bigg)\Bigg].~~~~~~~~~~~~~~~~~~~~~~~~~~
\end{eqnarray}
\end{example}

\begin{remark}\label{rmk: T-product extension}
The T-product tensor, akin to the Einstein product tensor, exhibits the spectrum decomposition property. Consequently, the findings presented in this study can be seamlessly extended to the T-product tensor, as demonstrated by the research conducted by~\cite{chang2021T-TensorExp}.
\end{remark}

\bibliographystyle{IEEETran}
\bibliography{TenExpOverManifold_Bib}

\end{document}